\documentclass[11pt]{amsart}
\usepackage{amsmath, amssymb, amsfonts, amsthm, amscd, latexsym}
\usepackage[mathscr]{eucal}

\theoremstyle{plain}
\newtheorem{lemma}{Lemma}[section]
\newtheorem{proposition}{Proposition}[section]
\newtheorem{theorem}{Theorem}[section]

\theoremstyle{definition}
\newtheorem{definition}{Definition}[section]

\newtheorem{remark}{Remark}[section]

\numberwithin{equation}{section}

\newcommand{\Aut}{{\rm Aut}}
\newcommand{\Co}{{\rm Co}}

\newcommand{\Hol}{{\rm Hol}}

\newcommand{\U}{{\rm U}}

\renewcommand{\O}{\Omega}
\renewcommand{\o}{\omega}

\newcommand{\grp}{{\mathbb G}}
\newcommand{\dgrp}{{\mathbb D}}

\newcommand{\Geod}{{\rm Geod}}
\newcommand{\geod}{{\rm geod}}

\newcommand{\obg}{{\rm Ob(\mathbb G)}}

\newcommand{\Osmooth}{{\Omega^{\infty}(X,*)}}

\newcommand{\gcalp}{{\mathbb G(\mathcal P)}}

\newcommand{\glob}{{\rm glob}}

\newcommand{\wti}{\widetilde}

\renewcommand{\a}{\alpha}
\newcommand{\be}{\beta}

\newcommand{\lra}{{\longrightarrow}}
\newcommand{\ra}{{\rightarrow}}

\begin{document}

\title[Connections, local subgroupoids, and a holonomy
Lie groupoid of a line bundle gerbe] {Connections, local
subgroupoids,\\ and a holonomy Lie groupoid \\ of a line bundle
gerbe}

\date{\today}

\author[R.~Brown and J.F.~Glazebrook]
{Ronald Brown and James F. Glazebrook}

\address{School of Mathematics\\
University of Wales\\ Dean St. Bangor\\ Gwynedd LL57 1UT UK}
\email[R. Brown]{r.brown@bangor.ac.uk}

\address{Department of Mathematics\\
Eastern Illinois University\\ Charleston IL 61920 USA \\ and
Department of Mathematics \\ University of Illinois at
Urbana--Champaign \\ Urbana IL 61801 USA}
\email[J. F. Glazebrook]{cfjfg@eiu.edu , glazebro@math.uiuc.edu}


\begin{abstract}

Our main aim is to associate a holonomy Lie groupoid to the
connective structure of an abelian gerbe. The construction has
analogies with a procedure for  the holonomy Lie groupoid of a
foliation, in using a locally Lie groupoid and a globalisation
procedure. We show that path connections and 2-holonomy on line
bundles may be formulated using the notion of a connection pair on
a double category, due to Brown-Spencer, but now formulated in
terms of double groupoids using the thin fundamental groupoids
introduced by Caetano-Mackaay-Picken. To obtain a locally Lie
groupoid to which globalisation applies, we use methods of local
subgroupoids as developed by Brown-$\dot{\rm I}$\c{c}en-Mucuk.
\end{abstract}

\maketitle


\section{Introduction}

On investigating the potential applications of double groupoids in
homotopy theory, Brown and Spencer in 1976 \cite{BS} developed the notion
of a \emph{connection pair} $(\Upsilon, \Hol)$ consisting of
\emph{the transport} $\Upsilon$ and \emph{holonomy} `$\Hol$',
which led to an equivalence of crossed modules with edge symmetric
double groupoids with special connections. The key `transport law'
for $\Upsilon$ used in this equivalence was abstracted from a law
for path connections on principal bundles due to Virsik \cite{Virsik},
applied to a connection pair on the double category $\Lambda T$ of Moore
paths on a topological category $T$. The relation of these ideas
with the connections of differential geometry has been
undeveloped. However, there is now a growing interest in
2-dimensional ideas in holonomy, particularly in those areas of
mathematics and mathematical physics  where the theory of
\emph{gerbes} plays a prominent role (see e.g. \cite{BD}
\cite{Barrett} \cite{Breen1} \cite{Bryl} \cite{MP}). From a
technical point of view, it is useful in the case where $T$ is a
Lie groupoid to move from the double category $\Lambda T$ as above
to a smooth double groupoid. Our first step is to use the
notion of the thin fundamental groupoid $\Lambda^1_1 (X)$ of a smooth
manifold $X$ (see e.g. \cite{CP} \cite{MP}) (\S 3). A major step is
to construct a smooth connection pair from the data of gerbes and
2-holonomy in abelian gerbes (\S 4,5).

\medbreak
To obtain a locally Lie groupoid we use in \S 6 the methods
of local subgroupoids and their holonomy Lie groupoids as in
\cite{BI1} \cite{BI2}.
In particular, conditions are given in \cite{BI2} for a connection pair
to yield  a local subgroupoid and so a locally Lie groupoid
$(\grp, W)$ with an associated  holonomy Lie groupoid
$\Hol(\grp,W)$, via a  \emph{Globalisation Theorem} of Aof--Brown
\cite{AOFB} \footnote{The origin of this theorem as in the work of
J. Pradines, is explained in \cite{AOFB}}. Brown--Mucuk \cite{BM}
showed that it recovered as a special case (and so with a new
universal property), the holonomy Lie groupoid of a foliation
where $\grp$ is the equivalence relation determined by the leaves
of the foliation. Two important points about this construction are
(i) $\Hol(\grp,W)$ comes with a universal property, and (ii) the
basic method of construction of $\Hol(\grp,W)$ involves an
algebraic framework for the intuition of `iteration of local
procedures', using Ehresmann's local smooth admissible sections.
It is of course this intuition which is behind holonomy
constructions in differential geometry.

\medbreak Thus in this note we establish:

\begin{theorem}\label{lsgtheorem}
We can associate to certain abelian gerbe data $(\mathcal P,
\mathcal A, \Geod)$ over a path connected manifold $X$, a local
subgroupoid $C(\mathcal P, \mathcal A, \Geod)$~. Relative to a strictly
regular atlas $\mathcal U(\mathcal P, \mathcal A, \Geod)$, there
exists a holonomy Lie groupoid $\Hol(\mathcal P, \mathcal A)$ with
base space $X$~.
\end{theorem}

The universal property satisfied by this holonomy groupoid will be
investigated in later work.


\section{Transport and holonomy in groupoids--some background}

\subsection{Connection and transport in a double category}

Firstly, let $T = (H, X,\a, \be, \circ, \epsilon)$ be a
topological category in which $\a, \be, : H \lra X$ are the
initial and final maps, respectively, $\circ$ denotes partial
composition, and $\epsilon : X \lra H$ is the unit function.
 A double category is specified by four related category
structures~:
\begin{equation}\label{dcat}
\begin{cases}
(D, H, \a_1, \be_1, \circ_1, \epsilon_1)~, &(D, V, \a_2, \be_2,
\circ_2, \epsilon_2 )~, \\ (V, X, \a_, \be, \circ, e)~, &(H, X,
\a_, \be, \circ, f)~,
\end{cases}
\end{equation}
of which each of the structures of the first row is compatible
with the other. For more details, see for instance \cite{BMosa}.
The elements of $D$ are called \emph{squares},
those of $H$ and $V$
\emph{the horizontal and vertical edges}, respectively, while $X$
consists of \emph{points}. A double category can be enhanced by the
abstract notion of a \emph{connection} as specified by a pair
$(\Upsilon, \Hol)$ in which $\Hol : V \lra H$ ({\it the holonomy})
is a functor of categories, and $\Upsilon : V \lra D$ (\emph{the
transport}) is a function, such that in the formalism of the
governing (higher dimensional) algebraic rules \cite{BS} (see also
\cite{Brown1} \cite{BI2} \cite{BMosa}), we have~:
\begin{itemize}
\item[(1)]
The bounding edges of $\Hol(a)$ and $\Upsilon(a)$, for $a : x \ra
y$ in $V$, are described by the diagram
\begin{equation}
\begin{CD}
x @>\Hol (a)>> y
\\ @Va VV   @VV e_y V
\\ y @> f_y >> y
\end{CD}
\end{equation}

\medbreak
\item[(2)]
\emph{The transport law} holds. That is, if $a, b \in V$ and $a
\circ b$ is defined, then
\begin{equation}\label{trans1}
\Upsilon(a \circ b) = \bmatrix \Upsilon(a) & \epsilon_2(\Hol(b))
\\ \\ \epsilon_1(b) & \Upsilon(b)
\endbmatrix
\end{equation}
\end{itemize}
A notable example of this construction (following \cite{BS}) is
\emph{the double category of Moore paths} $\Lambda T = (\Lambda H, H,
\Lambda X, X)$ on the topological category $T = (H,X)$~.
Here the squares are elements of $\Lambda H$, the horizontal edges
are $H$, the vertical edges are $\Lambda X$, and the set of points
is $X$~. Accordingly, a connection pair $(\Upsilon, \Hol)$ for
$\Lambda T$ consists of $(1)$, the transport $\Upsilon: \Lambda X
\lra \Lambda H$, and $(2)$, the holonomy $\Hol: \Lambda X \lra H$.
One aim is to realise \cite{BS} for double groupoids and connection
pairs $(\Upsilon, \Hol)$ in terms of the geometric data available.
In view of the growing interest in 2-dimensional structures in
differential geometry \cite{BD} \cite{MP}, we note the following:
\begin{theorem}\label{BStheorem}{\rm{ \cite{BMack1} \cite{BMosa}
\cite{BS}}}~
The following categories are equivalent~:
1) Crossed modules of groupoids, 2) $2$--Groupoids, and 3) edge symmetric double groupoids with special connection.
\end{theorem}
Structures 1) and 3) have been used extensively in homotopical
investigations (see the survey \cite{Brown1}) while crossed
modules of Lie groups are used in \cite{MP}, and are called Lie
2-groups.

\subsection{The groupoid $\grp(P,X)$ associated to a principal
$G$--bundle}

Let $X$ be a smooth connected manifold and consider a principal
$G$--bundle $\pi: P \lra X$, where $G$ is a Lie group. There is an
associated locally trivial groupoid over $X$ given by
\begin{equation}\label{grpd1}
 \grp(P,X) = P \times_G P \rightrightarrows X~,
\end{equation}
where for $u_1, u_2, \in P,~ g \in G$, we have $(u_1,u_2)g = (u_1
\cdot g, u_2 \cdot g)$, with equivalence classes satisfying the
multiplication rule $[u_1, u_2][u_3, u_4]=[u_1 \cdot g, u_4]$, for
which $u_3 = u_2 \cdot g$ in the fibre over $\pi(u_2) =
\pi(u_3)$~. Furthermore, if $u_0 \in P,~ x_0 = \pi(u_0)$, then
there are homeomorphisms and isomorphisms respectively, given by
\begin{equation}\label{grpd2}
\begin{aligned}
P &\lra  {P \times_G P} \vert_{x_0}~, ~~ u \mapsto [u, u_0]~, \\ G
&\lra  {P \times_G P} \vert^{x_0}_{x_0}~, ~~ g \mapsto [u_0, g
\cdot u_0]~,
\end{aligned}
\end{equation}
(see \cite{Mack1} Ch. II). The groupoid $\grp(P,X)$ (sometimes called the
\emph{Ehresmann symmetry groupoid}) will play a
significant role in all that follows.


\section{The thin path groupoid}

\subsection{Thin higher homotopy groups}

Our development here follows \cite{CP} \cite{MP} to which we refer
for further details. Let the set of all smooth $n$--loops in $X$ be
denoted by $\O_n^{\infty}(X,
*)$, where for $n=1$, we shall just write $\Osmooth$~. The product
$\ell_1 \star \ell_2$ of two $n$--loops and the inverse of an $n$--loop are
well--defined.
\begin{definition}
Two loops $\ell_1, \ell_2$ are called \emph{rank--$n$ homotopic} or
{\it thin homotopic}, denoted by $\ell_1 \overset{n}{\sim}
\ell_2$, if there exists an $\epsilon > 0$, and a homotopy $H:
[0,1] \times I^n \lra X$, satisfying~:
\begin{itemize}
\item[(1)]
$t_i \in [0, \epsilon) \cup (1-\epsilon,1] \Rightarrow~H(s, t_1,
\ldots, t_n) = *, ~ 1\leqslant i\leqslant n$~.

\medbreak
\item[(2)]
$s \in [0, \epsilon) \Rightarrow~H(s, t_1, \ldots, t_n) =
\ell_1(t_1, \ldots, t_n)$~.

\medbreak
\item[(3)]
$s \in (1-\epsilon,1] \Rightarrow~H(s, t_1, \ldots, t_n) =
\ell_2(t_1, \ldots, t_n)$~.

\medbreak
\item[(4)]
$H$ is smooth throughout its domain.

\medbreak
\item[(5)]
rank ~$DH_{(s,t_1, \ldots, t_n)} \leqslant n$, throughout its
domain.
\end{itemize}
\end{definition}
We denote by $\pi_n^n(X,*)$ the set of equivalences classes under
$\overset{n}{\sim}$ of  (thin) $n$--loops in $X$~. Observe that
$\pi_n^n(X,*)$ is abelian for $n \geqslant2$~. Also, for $\dim X
\leqslant n$, we have $\pi_n^n(X, *) = \pi_n(X, *)$, and for $\dim
X >n$, the group $\pi_n^n(X, *)$ is infinite dimensional.

\subsection{The smooth thin path groupoid $\Lambda^1_1(X)$}

Here we will set $n =1$~. For a smooth path $\lambda \in
\Lambda^1_1(X)$, a point $t_0 \in I$ is a \emph{sitting instant}
if there exists $\epsilon >0$ such that $\lambda$
is constant on $[0,1] \cap (t_0 - \epsilon, t_0 + \epsilon)$. It is
shown in \cite{CP} that there is always a re--parametrization of a
smooth path such that it sits this way at its endpoints. In this
case $\pi_1^1(X,*) = \Osmooth/{\overset{1}{\sim}}$~. Likewise,
there is the \emph{smooth thin path groupoid} $\Lambda^1_1(X)
\subset \Lambda(X)$ consisting of smooth paths $\lambda : I \lra X$
which are constant in a neighbourhood of $t=0,~t=1$, identified up
to rank $1$ homotopy, with
\begin{equation}
0 \leqslant t \leqslant \epsilon \Rightarrow H(s,t) = \lambda(0),~
1-
\epsilon \leqslant t \leqslant 1 \Rightarrow H(s,t) = \lambda(1)~,
\end{equation}
and with multiplication $\star$~. Henceforth, relevant path spaces
will be considered as smooth thin path groupoids.

\subsection{The transport law}

Following \cite{Virsik}, we outline several
properties of the (smooth) path connection
\begin{equation}
\begin{aligned}
\Upsilon~:~ \Lambda^1_1 (X) &\lra \Lambda^1_1 \grp(P,X)~, \\ \lambda
&\mapsto \Upsilon^{\lambda}~,
\end{aligned}
\end{equation}
which for $t \in [0,1]$, satisfies $\a  \Upsilon^{\lambda}(t) =
\lambda(0)$, and $\be  \Upsilon^{\lambda}(t) = \lambda(t)$~.

\medbreak If $\psi : [0,1] \lra [t_0, t_1] \subset [0,1]$ is a
diffeomorphism, we have the relationship
\begin{equation}\label{trans3}
\Upsilon^{\lambda} \cdot \psi =  \Upsilon^{\lambda \psi} \circ_2
\Upsilon^{\lambda}(\psi(0))~,
\end{equation}
which leads to $ \Upsilon^{\lambda}(0) = \wti{(\lambda(0))}$~.
Further, if $\lambda, \bar \lambda \in \Lambda^1_1(X)$ satisfy
$\lambda(1)= \bar \lambda(0)$, so that $\lambda \circ \bar \lambda$
is defined and is smooth, then
\begin{equation}
\lambda = (\lambda \circ \bar \lambda) \cdot \psi_1~,~ \bar
\lambda = (\lambda \circ \bar \lambda) \cdot \psi_2~,
\end{equation}
where $\psi_1(t) = \frac{1}{2}t$ and $\psi_2(t) = \frac{1}{2}t +
\frac{1}{2}$~. Given $\lambda \circ \bar \lambda \in \Lambda^1_1(X)$,
then on applying \eqref{trans3} to $\lambda \circ \bar \lambda$,
along with either $\psi_1$ or $\psi_2$, leads to an explicit
statement of the transport law for this case~:
\begin{equation}\label{trans4}
\Upsilon^{\lambda \circ \bar \lambda}(t) =
\begin{cases}& \Upsilon^{\lambda}(2t)~, ~0 \leqslant t \leqslant \frac{1}{2}~,
\\ & \Upsilon^{\bar \lambda}(2t-1) \circ_2
\Upsilon^{\lambda}(1)~,~ \frac{1}{2} \leqslant t \leqslant 1~.
\end{cases}
\end{equation}

In particular, we have $\Upsilon^{\lambda \circ \bar \lambda}(1) =
\Upsilon^{\bar \lambda}(1) \circ_2 \Upsilon^{\lambda}(1)$, and
$\Upsilon^{\lambda^{-1}}(1) = [\Upsilon^{\lambda}(1)]^{-1}$~.

\medbreak Now suppose that $\o_P$ denotes a given connection
$1$--form on $P \lra X$~. We refer to e.g. \cite{KN1} \cite{Mack1}
for the usual concept of horizontal path lifting and parallel
transport induced by $\o_P$ (which defines a `simple infinitesimal
connection' in the sense of \cite{Virsik}).
\begin{lemma}\label{pathlift}{\rm {\cite{Virsik} \cite{Mack1}
(Theorem $7.3$)}} Given a (smooth) principal $G$--bundle $P \lra
X$, a path connection $\Upsilon : \Lambda^1_1 (X) \lra \Lambda^1_1
\grp(P,X)$ is determined uniquely by a choice of connection
$1$--form $\o_P$~. Specifically, there exists a one--to--one
correspondence between
 $\Upsilon$ and $\o_P$, such that for $\gamma \in \Lambda^1_1(X)~,~\hat
\gamma = \Upsilon(\gamma)$ and $t_0 \in [0,1]$, we have
\begin{equation} \frac{d}{dt}~\hat \gamma(t_0) = (R_{\hat
\gamma(t_0)})_*~\o_P(\frac{d}{dt}\gamma(t_0))~.
\end{equation}
\end{lemma}
Thanks to thin homotopies, we can replace the topological double category
of Moore paths by a double Lie groupoid.
 Given a smooth principal $G$--bundle $P
\lra X$ with connection $1$--form $\o_P$, we can use
$\Lambda^1_1(X)$ together with the data provided by $\o_P$, to
produce a double groupoid~:
\begin{equation}\label{edgegrp1}
\begin{cases}
(\Lambda^1_1 \grp(P,X), \grp(P,X), \a_1, \be_1, \circ_1,
\epsilon_1)~, &(\Lambda^1_1 \grp(P,X), \Lambda^1_1(X), \a_2,
\be_2, \circ_2, \epsilon_2 )~, \\ (\Lambda^1_1(X), X, \a, \be,
\circ, e)~, & (\grp(P,X), X, \a, \be, \circ, e)~.
\end{cases}
\end{equation}
In the context of \cite{BS}, the existence of the connection pair
$(\Upsilon, \Hol)$ thus specializes to $(1)$, the parallel
transport as a smooth function on smooth groupoids $\Upsilon:
\Lambda^1_1(X) \lra \Lambda^1_1 \grp(P,X)$ satisfying the transport
law, and $(2)$, the holonomy $\Hol : \Lambda^1_1(X) \lra
\grp(P,X)$~.


\section{Gerbes and $2$--holonomy}

\subsection{Abelian gerbes}

The references for this section are \cite{Breen1} \cite{Bryl} \cite{Hitchin}
\cite{MP}. Let $X$ be a smooth (finite dimensional) connected
manifold and let $\mathcal U =\{U_i : i \in J\}$ be `good' open
cover of $X$ meaning that all $p$--fold ($p \geqslant 1$)
intersections $U_{i_1 \cdots i_p}
= U_{i_1} \cap \cdots \cap U_{i_p}$, are contractible.
The data for a {\it line bundle} ($G = \U(1)$) \emph{gerbe}
$\mathcal P
\lra X$, is given as follows~:
\begin{itemize}
\item
On each $U_{ij}$, there is a line bundle $L_{ij} \lra U_{ij}$,
satisfying $L_{ij} \cong L^{-1}_{ij}$~.

\medbreak
\item
There are trivializations $\theta_{ijk}$ of $L_{ij}L_{jk}L_{ki}$
on $U_{ijk}$ that satisfy the cocycle condition $\delta \theta
\equiv 1$ on $U_{ijkl}$~.
\end{itemize}
The corresponding data for a connective structure on $\mathcal P
\lra X$, is given as follows~:
\begin{itemize}
\item[(1)]
A $0$--\emph{object connection} is a covariant derivative
$\nabla_{ij}$ on $L_{ij}$, such that for each $i,j,k \in J$, it
satisfies the condition~:
\begin{equation}
\nabla_{ijk} \theta_{ijk} = 0~.
\end{equation}
In terms of the corresponding $1$--forms $A_{ij} \in A^1(U_{ij})$,
there is the equivalent relationship
\begin{equation}
\iota(A_{ij} + A_{jk} + A_{ki}) = - d \log g_{ijk}~,~ (A_{ij} = -
A_{ji})~,
\end{equation}
where $g$ is a \u{C}ech $2$--cocycle, $\delta g \equiv 1$~.

\medbreak
\item[(2)]
A $1$--\emph{connection} is defined by local $2$--forms $F_i \in
A^2(U_i)$ such that on $U_{ij}$, it satisfies
\begin{equation}
F_i - F_j = F_{ij} = \sigma^*(K(\nabla_{ij}))~,
\end{equation}
where $\sigma_{ij} \in \Gamma (U_{ij}, L_{ij})$, and $K$ denotes
the usual curvature. The latter is equivalent to the condition
$F_i - F_j = F_{ij} = dA_{ij}$~.
\end{itemize}
The abelian gerbe with its connective structure is denoted by
$(\mathcal P, \mathcal A)$~.

\subsection{The holonomy of $(\mathcal P, \mathcal A$)}

Suppose that $s : I^2 \lra X$ is a $2$--loop. The pull--back gerbe
$s^*(\mathcal P)$ is then a trivial gerbe and we can choose some
trivialization such that an object is given in terms of line
bundles trivialized by sections $\sigma_i$ over $V_i =
s^{-1}(U_i)$, with an object connection $\nabla_i$~. A global
$2$--form $\varepsilon$ (the error $2$--form) is defined on $I^2$
by
\begin{equation}
\varepsilon \vert V_i = s^*(F_i) - \sigma_i^*(K(\nabla_i))~.
\end{equation}
The holonomy $\Hol(s)$ of $(\mathcal P, \mathcal A)$ around the
$2$--loop $s$ is then given by
\begin{equation}
\exp (\iota ~\int_{I^2} \varepsilon)~. \end{equation} It is
independent of the choice of object and the connection on the
object, and is constant on thin homotopy classes. Furthermore,
$\Hol$ defines the $2$--\emph{holonomy} of $(\mathcal P,\mathcal
A)$ in terms of a  group homomorphism
\begin{equation}\label{Hol2}
\Hol ~:~ \pi_2^2(X,*) \lra \U(1)~,
\end{equation}
depending on $(\mathcal P, \mathcal A)$ up to equivalence, and
which is smooth on families of $2$--loops in $\O_2^{\infty}(X)$
when projected to $\pi^2_2(X,*)$ (see \cite{MP}).


\section{Parallel transport and holonomy in abelian gerbes}

\subsection{Transport of the gerbe data}

We have already noted the (thin) parallel transport
\begin{equation}
\Upsilon ~:~ \Lambda^1_1(X) \lra \Lambda^1_1 \grp(P,X )~,
\end{equation}
determined by the connection $1$--form $\o_P$ on a principal
$G$--bundle $P \lra X$~. Now we look for the analogous functor in
the case of the gerbe connection. The idea is that the gerbe data
$(\mathcal P, \mathcal A)$ determines a groupoid on the thin loop
groupoid $L_1^1(X) \subset \Lambda^1_1(X)$~. But if we assume $X$
is path connected and fix a base point, as we will henceforth,
then the gerbe data will readily lead us to the relevant transport
groupoids over their space of objects $\pi_1^1(X,*)$~.

\medbreak Using $(\mathcal P, \mathcal A)$,
 it is shown in \cite{Bryl} that there corresponds a smooth line bundle
$L^{\mathcal P} \lra \Osmooth$ with connection~. Here one
considers a quadruple of the type $(\gamma, F, \nabla, z)$ where
$\gamma \in \Osmooth$, $F$ is an object for $\gamma^*\mathcal P$
on $S^1$, $\nabla$ is an object connection in $F$, and $z \in
\mathbb C^*$~. These are defined up to a certain equivalence
relation. Now consider a homotopy $H : I^2 \lra X$ between
loops $\gamma, \mu$ and let $F, \nabla$ denote the object and
object connections respectively, for the pull--back gerbe
$H^*(\mathcal P) \lra I^2$. The parallel transport $\Upsilon$
along the homotopy $H$ is given explicitly by~:
\begin{equation}\label{PT1}
\Upsilon (H)(\gamma, \gamma^*F, \gamma^* \nabla, 1) = (\mu,
\mu^*F, \mu^* \nabla, 1) \exp(\iota \int_{I^2} \varepsilon)~.
\end{equation}
Note that the smooth line bundle with connection
\begin{equation}\label{bundle1}
p : (L^{\mathcal P}, D) \lra \Osmooth~,
\end{equation}
is representable as a principal $\U(1)$--bundle
\begin{equation}\label{bundle2}
\U(1) \hookrightarrow (L^{\mathcal P}, \o_P) \lra \Osmooth~,
\end{equation}
with connection $1$--form $\o_P$~. Parallel transport under $\o_P$
is defined along the cylindrical groupoid
\begin{equation}
\mathbb C(X,*) = \Lambda^1_1 \Osmooth \rightrightarrows \Osmooth~,
\end{equation}
where elements of $\mathbb C(X,*)$ are regarded as homotopies
between loops and whose morphisms are thin homotopy classes of
homotopies between loops (via based loops). At the same time
\eqref{bundle2} as determined by the gerbe data leads to the groupoid
\begin{equation}\label{gcalp}
\gcalp = \grp (\mathcal P, \mathcal A, X)= L^{\mathcal P}
\times_{\U(1)} L^{\mathcal P} \rightrightarrows \Osmooth~.
\end{equation}
Let then $T = (\gcalp, \Osmooth,\a, \be, \circ, e)$ whose
associated double category of paths $\Lambda^1_1 T$ contains the
horizontal and vertical groupoids $H = \gcalp \rightrightarrows
\Osmooth$, and $V = \mathbb C(X,*) \rightrightarrows \Osmooth$,
respectively. As a result $\Lambda^1_1 T$ is specified by the four
related groupoids~:
\begin{equation}
\begin{cases}
(\Lambda^1_1 \gcalp, \gcalp, \a_1, \be_1, \circ_1, \epsilon_1)~,~
& (\Lambda^1_1 \gcalp, \mathbb C(X,*), \a_2, \be_2, \circ_2,
\epsilon_2 )~, \\ (\mathbb C(X,*), \Osmooth, \a_, \be, \circ,
e)~,~ & (\gcalp, \Osmooth, \a_, \be, \circ, f)~.
\end{cases}
\end{equation}
To proceed, let $s_{\gamma},s_{\mu} \in \Gamma (\Osmooth,
L^{\mathcal P})$ be smooth sections and set $g_{\varepsilon} =
\exp(\iota \int_{I^2} \varepsilon) \in \mathbb C^{*}$, so that
\eqref{PT1} can be expressed as $\Upsilon (H)(s_{\gamma}) =
s_{\mu} \cdot g_{\varepsilon}$~. In this way, we can reduce matters
to considering the usual parallel transport in the $\U(1)$--bundle
$(L^{\mathcal P}, \o_P) \lra \Osmooth$ as induced by $\o_P$~.
Following Lemma \ref{pathlift}, $\o_P$ uniquely determines a smooth
(thin) parallel transport
\begin{equation}
\Upsilon : \mathbb C(X,*) \lra \Lambda^1_1\gcalp~,
\end{equation}
via the homotopy $H = \circ_2$ (horizontal structure) satisfying
the transport law \eqref{trans4}~. As for the holonomy, we see from
\eqref{grpd2} that the assignment $g \mapsto [u_0, u_0 \cdot g]$
for $u_0 \in L^{\mathcal P}, ~ p(u_0) = \gamma$, induces just as in
\eqref{grpd2} an isomorphism
\begin{equation}
\U(1) \cong \gcalp \vert_{\gamma}^{\gamma}~.
\end{equation}
Since $X$ is (path) connected, this leads to the holonomy  $\Hol :
\mathbb C(X,*) \lra \U(1)$, and in the context of the double
category of paths, the holonomy functor $\Hol : \mathbb C(X,*)
\lra \gcalp$~. It is straightforward to check that $\Hol(a)$
satisfies the relations showing that it is a bounding edge of the
square $\Upsilon(a)$~. We can summarize matters as follows~:
\begin{proposition}\label{gerbehol1}
Given the $\U(1)$--gerbe data $(\mathcal P,\mathcal A)$ over a
path--connected space $X$, there is an associated double groupoid
of thin paths
\begin{equation}
\begin{cases}
(\Lambda^1_1 \gcalp, \gcalp, \a_1, \be_1, \circ_1, \epsilon_1)~,~&
(\Lambda^1_1 \gcalp, \mathbb C(X,*), \a_2, \be_2, \circ_2,
\epsilon_2 )~,
\\ (\mathbb C(X,*), \Osmooth, \a_, \be, \circ, e)~,~&
(\gcalp, \Osmooth, \a_, \be, \circ, f)~,
\end{cases}
\end{equation}
and a connection pair $(\Upsilon,\Hol)$ given by $(1)$, the
transport $\Upsilon: \mathbb C(X,*) \lra \Lambda^1_1
\gcalp$, and $(2)$, the holonomy $\Hol :  \mathbb C(X,*) \lra \U(1)$~.
\end{proposition}

\subsection{Thin homotopies again}

If a pair of homotopies $H_1, H_2 : I^2 \lra X$, between a given
pair of paths, are themselves homotopic via a homotopy $Q : I^3
\lra X$, then the parallel transport around $H_1H_2^{-1}$ is
expressed by
\begin{equation}
\int_{I^3} Q^*B~.
\end{equation}
Accordingly, the parallel transport along $H_1$ and $H_2$ is the
same if the latter are thin homotopic since $Q$ may be chosen to
have rank $\leqslant 2$ everywhere. In this way we actually achieve
a line bundle descending to $\pi_1^1(X,*)$~:
\begin{equation}\label{brylinski}
p : (L^{\mathcal P}, D) \lra \pi_1^1(X,*)~.
\end{equation}
It will be convenient to express this in terms of the
principal $\U(1)$--bundle with connection $1$--form $\o_P$,
\begin{equation}\label{bundle3}
\U(1) \hookrightarrow (L^{\mathcal P}, \o_P) \lra \pi_1^1(X,*) ~,
\end{equation}
together with the groupoid
\begin{equation}\label{gerbegrp}
\gcalp = \grp(\mathcal P, \mathcal A, X) = L^{\mathcal P}
\times_{\U(1)} L^{\mathcal P} \rightrightarrows \pi_1^1(X,*)~.
\end{equation}

\medbreak Next we recall the cylindrical groupoid $\mathbb C(X,*)
\rightrightarrows \Osmooth$ and its morphisms  regarded as
thin homotopy classes of homotopies constituting the vertical
structure of $\Lambda^1_1 T$~. As explained in \cite{MP}, the
horizontal composition $\circ = \star$ determines the monoidal
composition of homotopies and  structure defined via the
composition of loops $\circ_1$, as well as the  corresponding
composition of vertical homotopies $H : I^2 \lra X$ between
concatenated loops $\lambda, \mu$, say. The descent to
$\pi_1^1(X,*)$ induces the \emph{thin cylindrical groupoid}
$\mathbb C_2^2(X,*) \rightrightarrows \pi_1^1(X,*)$~. Intuitively,
we can view the latter as given by
\begin{equation}
\mathbb C_2^2(X,*) = \Lambda^1_1(\pi_1^1(X,*)) =
\Lambda (\Osmooth/{\overset{1}{\sim}})/{\overset{1}{\sim}}~,
\end{equation}
which encapsulates the $2$--holonomy. In order to simplify the
notation, let us set $Y = \pi_1^1(X,*)$~. Consequently under the
relation ${\overset{1}\sim}$, the double category $\Lambda^1_1 T$
specializes to a double groupoid
\begin{equation}
\Lambda^1_1 T  = (\Lambda^1_1 \gcalp, \gcalp, \mathbb C_2^2(X,*),
Y)~,
\end{equation}
for which the squares are elements of $\Lambda^1_1 \gcalp$, the
horizontal edges are $\gcalp$, the vertical edges are $\mathbb
C_2^2(X,*)$, and the set of points is $Y = \pi_1^1(X,*)$~.
\begin{proposition}{\rm{(cf \cite{MP})}}
With respect to the principal $\U(1)$--bundle $(L^{\mathcal P},
\o_P) \lra \pi_1^1(X,*)$ determined by the gerbe data $(\mathcal
P, \mathcal A)$, we have a double groupoid
\begin{equation}
\begin{cases}
(\Lambda^1_1 \gcalp, \gcalp, \a_1, \be_1, \circ_1, \epsilon_1)~,~&
(\Lambda^1_1 \gcalp, \mathbb C_2^2(X,*), \a_2, \be_2, \circ_2,
\epsilon_2 )~,
\\ (\mathbb C_2^2(X,*), Y, \a_, \be, \circ, e)~,~
&(\gcalp, Y, \a_, \be, \circ, f)~,
\end{cases}
\end{equation}
and a connection pair $(\Upsilon,\Hol)$ given by $(1)$, the
transport $\Upsilon: \mathbb C_2^2(X,*) \lra \Lambda^1_1
\gcalp$, and $(2)$,
the holonomy $\Hol :  \mathbb C_2^2(X,*) \lra \gcalp$~. In
particular, $\gcalp$ and $\mathbb C_2^2(X,*)$ each admit the
structure of a Lie $2$--groupoid.
\end{proposition}
Effectively, this is a special case of Proposition \ref{gerbehol1}
when restricted to $Y=\pi_1^1(X,*)$~. The main point is the
existence of a certain normal monoidal subgroupoid $N(X,*)$ of
$\mathbb C(X,*)$ \cite{MP} such that on factoring--out by
${\overset{1}\sim}$, the horizontal arrows in the diagram below,
represent well--defined morphisms of groupoids~:
\begin{equation}
\begin{CD}
\mathbb C(X,*) @> {\overset{1}\sim}>> \mathbb C_2^2(X, *) =
\mathbb C(X,*)/{N(X,*)} \\
 @V{\downdownarrows}VV  @VV{\downdownarrows}V  \\
\Osmooth @> {\overset{1}\sim}>> Y = \Osmooth/{\overset{1}\sim}
 = \pi_1^1(X,*)~.
\end{CD}
\end{equation}
\begin{remark}
It is shown in \cite{MP} that the Lie $2$--groupoids
\begin{equation}\label{d1}
\dgrp_1 = (\mathbb C_2^2(X,*), H_1, V_1, Y, \circ, id)~, ~
\dgrp_2 = (\gcalp, H_2, V_2, Y, \star, id)~,~{\text{($V_1, V_2$
discrete)}}~,
\end{equation}
together with their respective monoidal structures, actually reduce
to Lie $2$--groups, and when $X$ is simply--connected,
the gerbe data $(\mathcal P,\mathcal A)$ can be constructed directly from the finer $2$--holonomy
$\Hol : \pi_2^2(X, *) \lra \U(1)$, and conversely.
\end{remark}


\section{Local subgroupoids}

\subsection{The local subgroupoid of a path connection}

In this section we describe how a holonomy Lie groupoid can be
associated to a $\U(1)$--gerbe using a local subgroupoid
constructed from its local connective structure. To proceed, let
$X$ be a topological space and $\a, \be : \grp \rightrightarrows X
= \obg$, a topological groupoid. For an open set $U \subset X$,
let $\grp \vert U$ be the full subgroupoid of $\grp$ on $U$~. Let
$L_{\grp}(U)$ be the set of all wide subgroupoids of $\grp \vert
U$~. For $V \subseteq U$, there is a restriction map $L_{UV} :
L_{\grp}(U) \lra L_{\grp}(V)$ sending $H \mapsto H \vert V$~. Thus
$L_{\grp}$ has the structure of a presheaf on $X$~.

\medbreak
Consider the sheaf $p_{\grp}~:  \mathcal L_{\grp} \lra
X$ formed from the presheaf $L_{\grp}$~. For $x \in X$, the stalk
$p_{\grp}^{-1}(x)$ of $\mathcal L_{\grp}$ has elements the germs
$[ U, H_U]_x$, for $U$ open in $X$ and $x \in U$~. $H\vert U$ is a
wide subgroupoid of $\grp \vert U$ and the equivalence relation
$\sim_x$ yielding the germ at $x$, is such that $H_U \sim_x K_V$,
where $K_V$ is a wide subgroupoid of $\grp \vert V$ if and only if
there exists a neighbourhood $W$ of $x$ such that $W \subseteq U
\cap V$ and $H_U \vert W = K_V \vert W$~.
 The topology of $\mathcal L_{\grp}$ is the usual sheaf topology
 with a sub--base of sets $\{ [U, H_U]_x : x \in X \}$, for all
 open subsets $U$ of $X$ and wide subgroupoids $H$ of $\grp \vert
 U$~.
\begin{definition}
 A \emph{local subgroupoid} of $\grp$ on the topological space
 $X$ is a continuous global section $S$ of the sheaf $p_{\grp}~:
 \mathcal L_{\grp} \lra X$ associated to the presheaf $L_{\grp}$~.
\end{definition}
Associated to a local subgroupoid are a number of technical
features such as the type of `atlas' with which one needs to work.
For instance, there is a `regular atlas' with the `globally
adapted' property; in the case of a  Lie local subgroupoid, a
`strictly regular atlas', etc. For an explanation of these terms
and further properties we refer to \cite{BI1} \cite{BIM}.

\medbreak
Suppose we have a continuous path connection $\Upsilon :
\Lambda^1_1(X) \lra \Lambda^1_1(\grp)$ with the usual properties as
before. We denote by $C_{\Upsilon}(\grp)$ the set of all $g \in
\grp$ such that if $\a(g) = x$, then there exists a path $\lambda$
in $X$ such that $\Upsilon(\lambda)$ joins $g$ to the identity
$\mathbf 1_x$; that is, $\Upsilon(\lambda)(0) = \mathbf 1_x$ and
$\Upsilon(\lambda)(1) = g$~. We next state some essential
properties of $C_{\Upsilon}(\grp)$ following \cite{BIM} \S 4~:
\begin{proposition}~
\begin{itemize}
\item[(1)]~$C_{\Upsilon}(\grp)$ is a wide subgroupoid of $\grp$~.

\medbreak
\item[(2)] If $\Upsilon$ is a path connection on $\grp$ and
$\mathcal U$ is an open cover of $X$, then $C_{\Upsilon}(\grp)$ is
generated by the family $C_{\Upsilon}(\grp\vert U)$, for all $U \in
\mathcal U$~.
\end{itemize}
\end{proposition}

\subsection{Geodesic and path local property of the
atlas}\label{geod}

In order to define a corresponding local subgroupoid
$C_{\Upsilon}(\grp, \mathcal U)$, it is necessary to work with an
atlas of the following type. Given an open cover $\mathcal U = \{
U_i : i \in I\}$ for $X$, we assume for each $i \in I$ there is a
collection of paths, denoted by $\geod(U_i)$ in $U_i$, whereby
$\lambda \in \geod(U_i)$ with $\lambda(0) =x, ~\lambda(1) = y$, is
called a \emph{geodesic path} from $x$ to $y$~. Further, we assume
\begin{itemize}
\item[(i)] If $x,y \in U_i$, then there is a unique geodesic path
$\geod_i(x,y)$ from $x$ to $y$~.

\medbreak
\item[(ii)] If $x,y \in U_i \cap U_j$, then $\geod_i(x,y)
=\geod_j(x,y)$~.

\medbreak
\item[(iii)] The path connection is flat for this structure,
meaning that if $\lambda : x \lra y$ is any path in $U_i$, then
$\Upsilon(\lambda)(1) = \Upsilon(\geod_i(x,y))(1)$~.
\end{itemize}
For such an atlas it follows from \cite{BIM} (Proposition $4.3$)
that there exists a local subgroupoid
\begin{equation}\label{lsg1}
C_{\Upsilon}(\grp, \mathcal U)(x) = [U_i, C_{\Upsilon}(\grp \vert
U_i)]_x~.
\end{equation}
We also need to specify the conditions to ensure that \eqref{lsg1}
can be globally adapted. Following \cite{BIM} (Proposition and
Definition 3.5 and 4.4), the equality  $C_{\Upsilon}(\grp,
\mathcal U)\vert U = C_{\Upsilon}(\grp \vert U, \mathcal U \cap
U)$ holds if for any $i,j \in I$ and $x \in U_i \cap U_j \cap U$,
there is an open set $W$ such that $x \in W \subseteq U_i \cap U_j
\cap U$, and $C_{\Upsilon}(\grp\vert U_i) \vert W =
C_{\Upsilon}(\grp \vert U_j \cap U) \vert W$~. Let us say that the
cover $\mathcal U$ is ($\Upsilon$--)\emph{path local for
$C_{\Upsilon}(\grp, \mathcal U)$} if this condition holds for all
open sets $U$ of $X$~. It follows from \cite{BIM} (Corollary 7.10)
that any ($\Upsilon$--) path local atlas of the local subgroupoid
$C_{\Upsilon}(\grp, \mathcal U)$, is globally adapted. Next let
\begin{equation}\label{sta}
W(\mathcal U_S) = \bigcup_{i \in I} H_i~,
\end{equation}
where we are given a strictly regular ($\Upsilon$--) path local
atlas $\mathcal U_S$ for $\grp$, and $H_i$ a Lie subgroupoid of
$\grp$~. Since such an atlas is globally adapted, it follows from
\cite{BI1} (Theorem 3.7) that there exists a locally Lie groupoid
$(\glob(C_{\Upsilon}(\grp\vert U)), W(\mathcal U_S))$~.
Furthermore, the Globalisation Theorem of \cite{BI1} (Theorem
$3.8$) establishes the existence of the associated holonomy Lie
groupoid $\Hol (S, \mathcal U_S)$~.

\subsection{Application to abelian gerbes--Proof of Theorem \ref{lsgtheorem}}

We proceed now to an application in the context of \cite{Bryl}
Chapter $5$ to which we refer for the notions of a \emph{torsor, a
connective structure on a sheaf of groupoids} as well as other
details. The application is also in the context of thin homotopies
as in \cite{MP}.

\medbreak
Let $\grp \rightrightarrows X$ be a groupoid and
$\mathcal U = (U_i)_{i \in I}$ be a good open covering of $X$~. As
above, we consider full subgroupoids $\grp \vert U_i$ as well as
wide subgroupoids $H_i$ of the latter. We consider principal
$G$--bundles $P_i \lra U_i$, along with isomorphisms
\begin{equation}
u_{ij} : P_j \vert U_{ij} ~\overset{\cong} \lra~ P_i \vert U_{ij}~,
\end{equation}
in $\grp \vert U_{ij}$~. As previously we assume that $G$ is the
abelian group $\U(1)$~. We consider a section $h_{ijk}$ (of the
band $\underline{\mathbb C}^*$) over $U_{ijk}$ by $h_{ijk} =
u_{ik}^{-1} u_{ij} u_{jk}$ where the latter is viewed as an
equality in $\Aut (P_k)$, noting that this corresponds to a
\u{C}ech $2$--cocycle. Let us decree the full subgroupoids
$\grp \vert U_i$ to be $\grp(P_i, U_i)$ and the $H_i$ to be locally
sectionable wide Lie subgroupoids of the latter.

\medbreak Just as before let $\mathcal L_{\grp} \lra X$ be the
sheaf corresponding to the presheaf $L_{\grp}$ of wide Lie
subgroupoids of $\grp$~. As a sheaf of groupoids in its own right,
we assume that $\mathcal L_{\grp}$ is equipped with a connective
structure $\Co$ (in the sense of \cite{Bryl}). Next we choose an
object $\o_i$ of the torsor $\Co(P_i)$, where we regard $\o_i$ as
simply a connection $1$--form on $P_i \lra U_i$, and assume the
geodesic--path local property (\S \ref{geod}) of $\mathcal U$
relative to the $\o_i$~. We denote by $(\mathcal P, \mathcal A, \Geod)$
 the corresponding abelian gerbe data together with this property~. The next step is to apply the techniques of \S 3 to this situation~.

\medbreak
To proceed, we define a $1$--form $\o_{ij}$ on $U_{ij}$
by
\begin{equation}
\o_{ij} = \o_i - (u_{ij})_{*}(\o_j)~.
\end{equation}
Recalling $h_{ijk} = u_{ik}^{-1} u_{ij} u_{jk}$, it follows that
\begin{equation}
\o_{ij} + \o_{jk} - \o_{ik} = \o_i - (u_{ij} u_{jk}
u_{ki})_{*}(\o_i) = h_{ijk}^{-1}~dh_{ijk} ~.
\end{equation}
This data, denoted $(\underline{h}, \underline{\o})$, so defines a
\u{C}ech $2$--cocycle, but with coefficients in the complex of
sheaves $\underline{\mathbb C}^* ~\overset{d\log} \lra
~\underline{A}^1_{X, \mathbb C}$~. On restriction to thin path
groupoids, the object connections $\o_i$ of $\Co(P_i)$ determine on each $U_i$
a (thin) path connection
\begin{equation}
\Lambda^1_1(U_i) \overset{\Upsilon_i} \lra \Lambda^1_1(H_i) \subset
\Lambda^1_1 (\grp\vert U_i)~,
\end{equation}
satisfying the local flatness property $\Upsilon_i(\lambda)(1) =
\Upsilon_i (\geod_i(x,y))(1)$~. Also, there is an open set $W
\subseteq U_{ij}$, for which $H_i \vert W = H_j \vert W = H_{ij}$,
and so on the overlaps $U_{ij}$, we have a (thin) path connection
\begin{equation}
\Lambda^1_1(U_{ij}) \overset{\Upsilon_{ij}} \lra \Lambda^1_1(H_{ij})
\subset \Lambda^1_1(\grp\vert U_{ij})~.
\end{equation}
Consider the local subgroupoid of the atlas as given by $S(x) =
[U_i, H_i]_x$~. At the same time there is a wide subgroupoid
generated by the family $C_{\Upsilon_i}(\grp
\vert U_i)$ for all $U_i \in \mathcal U$~. As noted earlier, this
leads to a local subgroupoid $C(\mathcal P, \mathcal A, \Geod)$ given by
\begin{equation}
C(\mathcal P, \mathcal A, \Geod)(x)
= [U_i, C_{\Upsilon_i}(\grp\vert U_i)]_x~.
\end{equation}
Relative to a strictly regular path local atlas $\mathcal
U(\mathcal P, \mathcal A,\Geod)$, we apply the same considerations
as before along with the globalisation \cite{BI1} (Theorem $3.8$),
to obtain a holonomy Lie groupoid $\Hol(\mathcal P,
\mathcal A)$, thus establishing Theorem \ref{lsgtheorem}.


\end{document}